\documentclass{article}
\usepackage{eurosym}
\usepackage{amsmath}
\usepackage{amsfonts}
\usepackage{amssymb}
\usepackage{graphicx}
\usepackage{cite}

\newtheorem{theorem}{Theorem}[section]

\newtheorem{case}[theorem]{Case}

\newtheorem{definition}[theorem]{Definition}

\input{tcilatex}
\begin{document}

\title{On the Comparison of Perturbation-Iteration Algorithm and Residual
Power Series Method to Solve Fractional Zakharov-Kuznetsov Equation}
\author{Mehmet \c{S}ENOL and \ Hamed Daei Kasmaei \\
Nev\c{s}ehir Hac\i\ Bekta\c{s} Veli University, Department of Mathematics,
Nev\c{s}ehir, Turkey\\
Department of Mathematics and Statistics, Faculty of Science,\\
Central Tehran Branch, Islamic Azad University, Tehran, Iran.\\
e-mail: msenol@nevsehir.edu.tr, hamedelectroj@gmail.com}
\maketitle

\begin{abstract}
In this paper, we present analytic-approximate solution of a fractional
Zakharov-Kuznetsov equation by means of perturbation-iteration algorithm
(PIA) and residual power series method (RSPM). Basic definitions of
fractional derivatives are described in the Caputo sense. Several examples
are given and the results are compared to exact solutions. The results show
that both methods are competitive, effective, convenient and simple to use.

\end{abstract}

\section{Introduction}

Fractional differential equations have received considerable interest in
recent years and have been extensively investigated and applied for many
real problems which are modeled in various areas and have been the focus of
many studies due to their frequent appearance in various applications such
as fluid mechanics, viscoelasticity, biology, physics and engineering.
Therefore, great attention has been given to find numerical and
Analytic-approximate solutions of FDEs. Some of the recent
analytical-approximate methods for FDEs include the Adomian decomposition
method (ADM), the homotopy perturbation method (HPM), the variational
iteration method (VIM) and homotopy analysis method (HAM). The ADM was
applied to fractional diffusion equations in \cite{1} and fractional
modified KdV equations in \cite{2}. Hosseinnia et al. \cite{3} suggested an
enhanced HPM for FDEs and also Abdulaziz et al. \cite{4} improved the
application of HPM to systems of FDEs. In \cite{5}, Abdulaziz et al. solved
the fractional IVPs by the HPM. The HAM was applied to fractional
KDV-Burgers-Kuromoto equations \cite{6}, fractional IVPs \cite{7},
time-fractional PDEs \cite{8}, linear and nonlinear FDEs \cite{9} and
systems of nonlinear FDEs \cite{10}. Variational iteration method was
applied to solve some types of FDEs in \cite{12,13}.\newline
In this paper, we first introduce fractional Zakharov-Kuznetsov equation,
then we describe perturbation-iteration algorithm and residual power series
method (RSPM) in order to implement on Zakharov-Kuznetsov equation, then we
present some examples that show reliability and efficiency of two methods in
order to compare their numerical results. At last, we discuss about obtained
results as a section for conclusion. This paper considers the fractional
version of the Zakharov-Kuznetsov equations as studied in \cite{18}. The
fractional Zakharov-Kuznetsov equations considered are of the form:

\begin{equation}
D_{t}^{\alpha }u+a(u^{p})_{x}+b(u^{q})_{xxx}+b(u^{r})_{ttx}=0  \label{1}
\end{equation}%
where $u=u(x,y,t),$ $\alpha $ is a parameter describing the order of the
fractional derivative ($0<\alpha \leq 1$), $a$ and $b$ are arbitrary
constants and $p$, $q$ and $r$ are integers and $p;q;r\neq 0$ governs the
behavior of weakly nonlinear ion acoustic waves in a plasma comprising cold
ions and hot isothermal electrons in the presence of a uniform magnetic
field \cite{19,20}. The Zakharov-Kuznetsov equation was first derived for
describing weakly nonlinear ion-acoustic waves in strongly magnetized
lossless plasma in two dimensions \cite{21}. The FZK equations have been
studied previously by using VIM \cite{22} and HPM \cite{23}.

\section{Basic Definitions}

\begin{definition}
\bigskip A real function $f(t)$, $t>0$ is said to be in the space $C_{\mu }$%
, $(\mu >0)$ if there exists a real number $p(>\mu $) such that $%
f(t)=t^{p}f_{1}(t)$ where $f_{1}\in C[0,\infty )$ and it is said to be in
the space $C_{\mu }^{m}$ if $f^{(m)}\in C_{\mu }, m\in \aleph$ \cite{24}.
\end{definition}

\begin{definition}
The Riemann-Liouville fractional integral operator $(J^{\alpha })$ of order $%
\alpha \geq 0$ of a function $f\in C_{\mu }$, $\mu \geq -1$ is defined as 
\cite{25}:%
\begin{equation}
J^{\alpha }f(t)=\frac{1}{\Gamma (\alpha )}\dint\limits_{0}^{t}(t-\tau
)^{\alpha -1}f(\tau )d\tau ,\text{ \ \ }\alpha ,t>0  \label{2}
\end{equation}%
and $J^{0}f(t)=f(t)$, where $\Gamma $ is the well-known gamma function. For $%
f\in C_{\mu }$, $\mu \geq -1$, $\alpha $,$\beta \geq 0$ and $\lambda >-1$,
the following properties hold.

\begin{itemize}
\item $J^{\alpha }J^{\beta }f(t)=J^{\alpha +\beta }f(t)$

\item $J^{\alpha }J^{\beta }f(t)=J^{\beta }J^{\alpha }f(t)$

\item $J^{\alpha }t^{\lambda }=\frac{\Gamma (\lambda +1)}{\Gamma (\lambda
+1+\alpha )}t^{\lambda+\alpha }$.
\end{itemize}
\end{definition}

\begin{definition}
The Caputo fractional derivative of $f$ of order $\alpha $, $f\in C_{-1}^{m}$%
, $m\in \aleph \cup \{0\}$, is defined as \cite{26}:%
\begin{equation}
D^{\alpha }f(t)=J^{m-\alpha }f^{(m)}(t)=\frac{1}{\Gamma (m-\alpha )}%
\dint\limits_{0}^{t}(t-\tau )^{m-\alpha -1}f^{(m)}(\tau )d\tau ,\text{ \ \ }%
\alpha ,t>0  \label{3}
\end{equation}%
where $m-1<\alpha <m$ with the following properties;

\begin{itemize}
\item $D^{\alpha }\left( af(t)+bg(t)\right) =aD^{\alpha }f(t)+bD^{\alpha
}g(t),$ $a,b\in \Re $,

\item $D^{\alpha }J^{\alpha }f(t)=f(t)$,

\item $J^{\alpha }D^{\alpha }f(t)=f(t)-\dsum\limits_{j=0}^{k-1}f^{(j)}(0)%
\frac{t^{j}}{j!}$, $t>0$.
\end{itemize}
\end{definition}

\section{Overview of the Perturbation-Iteration Algorithm PIA(1,1)}

As one of the most practical subjects of physics and mathematics,
differential equations create models for a number of problems in science and
engineering to give an explanation for a better understanding of the events.
Perturbation methods have been used for this purpose for over a century \cite%
{27,28,29}.

But the main difficulty in the application of perturbation methods is the
requirement of a small parameter or to install a small artificial parameter
in the equation. For this reason, the obtained solutions are restricted by a
validity range of physical parameters. Therefore, to overcome the
disadvantages come with the perturbation techniques, some methods have been
suggested by several authors \cite{30,31,32,33,34,35,36,37,38}.

Parallel to these studies, recently a new perturbation-iteration algorithm
has been proposed by Aksoy, Pakdemirli and their co-workers \cite{39,40}. In
the new technique, an iterative algorithm is established on the perturbation
expansion.First the method applied for Bratu type second order equations 
\cite{39} to obtain approximate solutions. Then the algorithms were tested
on some nonlinear heat equations also \cite{40}. The solutions of the
Volterra and Fredholm type integral equations \cite{42}, ordinary
differential equation and systems \cite{43} and the solutions of ordinary
fractional differential equations \cite{44} have given by the present
method. Modification of the PIA has been also introduced by Bildik and Deniz 
\cite{49,50,51}.

In this study, the previously developed technique is applied to systems of
nonlinear fractional differential equations for the first time. To obtain
the approximate solutions of equations, the most basic
perturbation-iteration algorithm PIA(1,1) is employed by taking one
correction term in the perturbation expansion and correction terms of only
first derivatives in the Taylor series expansion, i.e. $n=1,$ $m=1$.

Consider the following initial value problem.

\begin{equation}
F\left( u_{t}^{(\alpha )},u_{x},u_{xxx},u_{ttx},\varepsilon \right) =0
\label{4}
\end{equation}%
\begin{equation*}
u(x,y,0)=c
\end{equation*}%
where $u=u(x,y,t)$ and $\varepsilon $ is a small perturbation parameter. The
perturbation expansions with only one correction term is

\begin{equation}
u_{n+1}=u_{n}+\varepsilon (u_{c})_{n}  \label{5}
\end{equation}%
where subscript $n$ represents the $n-th$ iteration.

Replacing Eq.(5) into Eq.(4) and writing in the Taylor Series expansion for
first order derivatives in the neighborhood of $\varepsilon =0$ gives 
\begin{eqnarray}
&&F\left( (u_{n}^{(\alpha
)})_{t},(u_{n})_{x},(u_{n})_{xxx},(u_{n})_{ttx},0\right) +F_{u_{t}^{(\alpha
)}}\left( (u_{n}^{(\alpha
)})_{t},(u_{n})_{x},(u_{n})_{xxx},(u_{n})_{ttx},0\right) \varepsilon \left(
(u_{c}^{(\alpha )})_{t}\right) _{n}  \notag \\
&&+F_{u_{x}}\left( (u_{n}^{(\alpha
)})_{t},(u_{n})_{x},(u_{n})_{xxx},(u_{n})_{ttx},0\right) \varepsilon \left(
(u_{c})_{x}\right) _{n}+F_{u_{xxx}}\left( (u_{n}^{(\alpha
)})_{t},(u_{n})_{x},(u_{n})_{xxx},(u_{n})_{ttx},0\right) \varepsilon \left(
(u_{c})_{xxx}\right) _{n}  \notag \\
&&+F_{u_{ttx}}\left( (u_{n}^{(\alpha
)})_{t},(u_{n})_{x},(u_{n})_{xxx},(u_{n})_{ttx},0\right) \varepsilon \left(
(u_{c})_{ttx}\right) _{n}+F_{\varepsilon }\left( (u_{n}^{(\alpha
)})_{t},(u_{n})_{x},(u_{n})_{xxx},(u_{n})_{ttx},0\right) \varepsilon =0
\label{6}
\end{eqnarray}%
or%
\begin{equation}
\frac{F}{\varepsilon }+\left( (u_{c}^{(\alpha )})_{t}\right)
_{n}F_{u_{t}^{(\alpha )}}+\left( (u_{c})_{x}\right) _{n}F_{u_{x}}+\left(
(u_{c})_{xxx}\right) _{n}F_{u_{xxx}}+\left( (u_{c})_{ttx}\right)
_{n}F_{u_{ttx}}+F_{\varepsilon }=0  \label{7}
\end{equation}%
where $F_{u_{t}^{(\alpha )}}=\frac{\partial F}{\partial u_{t}^{(\alpha )}},$ 
$F_{u_{x}}=\frac{\partial F}{\partial u_{x}},$ $F_{u_{xxx}}=\frac{\partial F%
}{\partial u_{xxx}},$ $F_{u_{ttx}}=\frac{\partial F}{\partial u_{xxx}}$ and $%
F_{\varepsilon }=\frac{\partial F}{\partial \varepsilon }.$

All derivatives in the expansion are calculated at $\varepsilon =0$.
Therefore in the procedure of computations, each term is obtained when $%
\varepsilon $ tends to zero. The method converges in few iterations and in
fact we have a saturated solution after doing computations even in the
initial steps to find favorite approximate solution. Beginning with an
initial function $u_{0}(x,y,t)$, first $(u_{c})_{0}(x,y,t)$ has been
determined by the help of Eq.(7). Then using Eq.(5), $\left( n+1\right) .$
iteration solution could be found. Iteration process is repeated using
Eq.(7) and Eq.(5) until achieving an acceptable result. The ability of the
method is so high that it can become convergent in just a few of
computational iterations. The reliability and effectiveness of the method is
shown by an example after presenting the convergence of PIA method. The
purpose of this paper is to obtain approximate solutions of the fractional
Zakharov-Kuznetsov equations by RSPM and PIA and to determine series
solutions with high accuracy.

\section{Algorithm of RPSM}

In this section, we employ our technique of the RPS method to find out
series solution for the Zhakarov-Kunetsov equation. The RPS method \cite%
{45,46} consists of expressing the solution of (1) as a fractional
power-series expansion about the initial point $t=0$.\newline
\begin{equation}
u(x,y,t)=\sum\limits_{n=0}^{\infty }{{{f}_{n}}(x,y)\frac{{{t}^{n\alpha }}}{%
\Gamma \left( 1+n\alpha \right) }}\text{ },\text{ }0<\alpha \leq 1,\text{ }%
x\in \left[ a,b\right] ,\text{ }0\leq t<\Re   \label{8}
\end{equation}%
Next, we let $u_{k}(x,t)$ to denote the $k$-th truncated series of $u(x,t)$,
i.e., 
\begin{equation}
u_{k}(x,y,t)=\sum\limits_{n=0}^{k}{{{f}_{n}}(x,y)\frac{{{t}^{n\alpha }}}{%
\Gamma \left( 1+n\alpha \right) }}\text{ },\text{ }0<\alpha \leq 1,\text{ }%
x\in \left[ a,b\right] ,\text{ }0\leq t<\Re   \label{9}
\end{equation}%
To achieve our goal, we suppose that the zeroth RPS approximate solutions of 
$u(x,t)$ is as follows : 
\begin{equation}
{{u}_{0}}(x,y,t)={{f}_{0}}(x,y)=u(x,y,0)=f(x,y)  \label{10}
\end{equation}%
Also, Eq.(13) can be written as : 
\begin{equation}
{{u}_{k}}\left( x,y,t\right) =f(x,y)+\sum\limits_{n=1}^{k}{{{f}_{n}}(x,y)%
\frac{{{t}^{n\alpha }}}{\Gamma \left( 1+n\alpha \right) }},\text{ }0<\alpha
\leq 1\text{, }x\in \left[ a,b\right] ,\text{ }0\leq t<R,\text{ }%
k=1,2,3,\cdots   \label{11}
\end{equation}%
Now, we define the residual functions, $Res$, for Eq.(1) as 
\begin{equation}
{Res}_{u}\left( x,y,t\right) =D_{t}^{\alpha }u+a{{\left( {{u}^{p}}\right) }%
_{x}}+b{{\left( {{u}^{q}}\right) }_{xxx}}+b{{\left( {{u}^{r}}\right) }_{yyx}}
\label{12}
\end{equation}%
and, therefore, the k-th residual function, $Res_{u,k}(x,t)$ is given as an
iterative relation as 
\begin{equation}
{{Res}_{u,k}}\left( x,y,t\right) =D_{t}^{\alpha }{{u}_{k}}+a{{\left(
u_{k}^{p}\right) }_{x}}+b{{\left( u_{k}^{q}\right) }_{xxx}}+b{{\left(
u_{k}^{r}\right) }_{yyx}}  \label{13}
\end{equation}%
substitution of $k$-th truncated series $u_{k}(x,y,t)$ of Eq.(13) into
Eq.(16) leads to the following definition for the $k$-th residual function
that is called $Res^{k}(x,y,t)$ and in this case regarding the main Eq.(1),
we have :

\begin{align}
& {{R}^{k}}(x,y,t)=D_{t}^{\alpha }\left( f(x,y)+\sum\limits_{n=1}^{k}{{{f}%
_{n}}(x,y)\frac{{{t}^{n\alpha }}}{\Gamma \left( 1+n\alpha \right) }}\right)
+a\left( f(x,y)+\sum\limits_{n=1}^{k}{{{f}_{n}}(x,y)\frac{{{t}^{n\alpha }}}{%
\Gamma \left( 1+n\alpha \right) }}\right) _{x}^{p}  \notag \\
& +b\left( f(x,y)+\sum\limits_{n=1}^{k}{{{f}_{n}}(x)\frac{{{t}^{n\alpha }}}{%
\Gamma \left( 1+n\alpha \right) }}\right) _{_{xxx}}^{q}+b\left(
f(x,y)+\sum\limits_{n=1}^{k}{{{f}_{n}}(x,y)\frac{{{t}^{n\alpha }}}{\Gamma
\left( 1+n\alpha \right) }}\right) _{ttx}^{r}  \label{14}
\end{align}

and we have : 
\begin{equation}
Res^{\infty }(x,y,t)=\underset{k\rightarrow \infty }{\mathop{\lim }}\,{{Res}%
^{k}}(x,y,t)  \label{15}
\end{equation}

It is easy to see that $Res^{\infty }(x,y,t)=0$ for each $x\in \lbrack a,b]$.
This show that $Res^{\infty }(x,y,t)$ is infinitely many times
differentiable at $x=a$. On the other hand, $\frac{{{d}^{k}}}{d{{x}^{k-1}}}{{%
R}^{\infty }}(x,y,0)=\frac{{{d}^{k}}}{d{{x}^{k-1}}}{{R}^{k}}(x,y,0)=0$. In
fact, this relation is a fundamental rule in RPS method and its applications
[13].\newline
Now, in order to derive the $k$-th approximate solution, we consider $%
u_{k}(x,y,t)=\sum\limits_{n=0}^{k}{{{f}_{n}}(x,y)\frac{{{t}^{n\alpha }}}{%
\Gamma \left( 1+n\alpha \right) }}$ and we differentiate both sides of
Eq.(13) with respect to $x,$ $y,$ $t$\ and substitute $t=0$ in order to find
constant parameters. After substituting these parameters in $u_{k}(x,y,t)$,
we can obtain $k$-th truncated series and by putting it in Eq.(1), we reach
to our favorite approximate solution. This procedure can be repeated till
the arbitrary order coefficients of RPS solutions for Eq.(1). Moreover,
higher accuracy can be achieved by evaluating more components of the
solution.

\section{Applications}

\begin{case}
PIA: Consider the following time-fractional FZK equation \cite{47,48}:
\end{case}

\begin{equation}
D_{t}^{\alpha }u(x,y,t)+(u^{2}(x,y,t))_{x}+\frac{1}{8}(u^{2}(x,y,t))_{xxx}+%
\frac{1}{8}(u^{2}(x,y,t))_{xyy},\ \ \ t>0,\ \ \ 0\leq t<1,\ \ \ 0<\alpha
\leq 1  \label{16}
\end{equation}%
with the initial condition\ $u\left( x,y,0\right) =\frac{4}{3}\rho \sinh
^{2}(x+y)$ and the known exact solution for $\alpha =1$ is

\begin{equation}
u\left( x,y,t\right) =\frac{4}{3}\rho \sinh ^{2}(x+y-\rho t)  \label{17}
\end{equation}

Before iteration process rewriting Eq.(16) with adding and subtracting $%
u_{t}(x,y,t)$ and inserting artificial parameter $\varepsilon $ to the
equation gives%
\begin{eqnarray}
F\left( u^{^{{\prime }}},u,\varepsilon \right)  &=&\varepsilon \frac{%
d^{\alpha }u_{n}(x,y,t)}{{dt}^{\alpha }}+\frac{\partial }{\partial t}%
u_{n}(x,y,t)-{\varepsilon }\frac{\partial }{\partial t}%
u_{n}(x,y,t)+2u_{n}(x,y,t)\frac{\partial }{\partial x}u_{n}(x,y,t)  \notag \\
&&+\frac{3}{4}\frac{\partial }{\partial x}u_{n}(x,y,t)\frac{\partial ^{2}}{%
\partial x^{2}}u_{n}(x,y,t)+\frac{1}{4}u_{n}(x,y,t)\frac{\partial ^{3}}{%
\partial x^{3}}u_{n}(x,y,t)  \notag \\
&&+\frac{1}{4}\frac{\partial ^{2}}{\partial y^{2}}u_{n}(x,y,t)\frac{\partial 
}{\partial x}u_{n}(x,y,t)+\frac{1}{2}\frac{\partial }{\partial y}u_{n}(x,y,t)%
\frac{\partial ^{2}}{\partial x\partial y}u_{n}(x,y,t)  \notag \\
&&+\frac{1}{4}u_{n}(x,y,t)\frac{\partial ^{3}}{\partial x\partial y\partial y%
}u_{n}(x,y,t)  \label{18}
\end{eqnarray}%
and for the iteration formula%
\begin{equation}
u_{t}(x,y,t)+\frac{F_{u}}{F_{u_{t}}}u\left( x,y,t\right) =-\frac{%
F_{\varepsilon }+\frac{F}{\varepsilon }}{F_{u_{t}}}  \label{19}
\end{equation}%
the terms that will be replaced in, are

\begin{eqnarray}
F &=&(u_{n})_{t}(x,y,t)-1+\frac{t}{4}  \notag \\
F_{u} &=&0  \notag \\
F_{u_{t}} &=&1  \notag \\
F_{\varepsilon } &=&-{u_{n}^{^{{\prime }}}\left( x,y,t\right) }+\frac{1}{%
\Gamma (1-\alpha )}\varepsilon \int_{0}^{t}{\frac{(u_{n})_{t}(x,y,s)}{{(t-s)}%
^{\alpha }}ds}+2u_{n}(x,y,t)(u_{n})_{x}(x,y,t)+\frac{3}{4}%
(u_{n})_{x}(x,y,t)(u_{n})_{xx}(x,y,t)  \notag \\
&&+\frac{1}{4}u_{n}(x,y,t)(u_{n})_{xxx}(x,y,t)+\frac{1}{4}%
(u_{n})_{yy}(x,y,t)(u_{n})_{x}(x,y,t)+\frac{1}{2}%
(u_{n})_{y}(x,y,t)(u_{n})_{xy}(x,y,t)  \notag \\
&&+\frac{1}{4}u_{n}(x,y,t)(u_{n})_{xyy}(x,y,t)  \label{20}
\end{eqnarray}

After substitution the differential equation for this problem, Eq.(16)
becomes \bigskip

\begin{eqnarray}
&&\frac{4}{\Gamma (1-\alpha )}\int_{0}^{t}{{\left( t-s\right) }^{-\alpha
}(u_{n})_{t}(x,y,s)ds+}4u_{t}(x,y,t)+2(u_{n})_{y}(x,y,t)(u_{n})_{xy}(x,y,t) 
\notag \\
&&+(u_{n})_{x}(x,y,t)\left( (u_{n})_{yy}(x,y,t)+3(u_{n})_{xx}(x,y,t)\right) 
\notag \\
&&+u_{n}(x,y,t)\left(
8(u_{n})_{x}(x,y,t)+(u_{n})_{xyy}(x,y,t)+(u_{n})_{xxx}(x,y,t)\right) =0
\label{21}
\end{eqnarray}

Appropriate to the initial condition, chosen $u\left( x,y,0\right) =\frac{4}{%
3}\rho \sinh ^{2}(x+y)$ and solving Eq.(21) for $n=0$ gives%
\begin{equation}
{{u}_{c}((x,y,t))}_{0}=\frac{8}{9}t\left( 4\rho ^{2}\sinh (2x+2y)-5\rho
^{2}\sinh (4x+4y)\right) +c_{1}(x,y)  \label{22}
\end{equation}

This expression written in iteration formula

\begin{equation}
u_{1}=u_{0}+\varepsilon {{u}_{c}((x,y,t))}_{0}  \label{23}
\end{equation}%
yields

\begin{equation}
u_{1}\left( x,y,t\right) =u_{0}\left( x,y,t\right) +\varepsilon \left( \frac{%
8}{9}t\left( 4\rho ^{2}\sinh (2x+2y)-5\rho ^{2}\sinh (4x+4y)\right)
+c_{1}(x,y)\right)   \label{24}
\end{equation}%
or

\begin{equation}
u_{1}\left( x,y,t\right) =\varepsilon (t-\frac{t^{2}}{8}+C_{1})  \label{25}
\end{equation}

Solving this equation for initial condition

\begin{equation}
u_{1}\left( x,y,0\right) =\frac{4}{3}\rho \sinh ^{2}(x+y)  \label{26}
\end{equation}%
we obtain

\begin{equation}
c_{1}(x,y)=0  \label{27}
\end{equation}

For this value and $\varepsilon =1$ reorganizing $u_{1}(x,y,t)$

\begin{equation}
u_{1}(x,y,t)=-\frac{4}{9}\rho \left( 4t\rho \left( \cosh (x+y)+5\cosh
(3(x+y)\right) )-3\sinh (x+y))\right) \sinh (x+y)  \label{28}
\end{equation}%
gives the first iteration result. If the iteration procedure is continued in
a similar way, we obtain the following\ second iteration.

\begin{eqnarray}
u_{2}(x,y,t) &=&\frac{2}{243}\rho (\frac{108\rho t^{2-\alpha }(5\sinh
(4(x+y))-4\sinh (2(x+y)))}{\gamma }{(3-\alpha )}  \notag \\
&&+9(208\rho ^{2}t^{2}+9)\cosh (2(x+y))-8\rho t(10\rho t(8\rho t(4\sinh
(2(x+y))+8\sinh (4(x+y))  \notag \\
&&-60\sinh (6(x+y))+85\sinh (8(x+y)))+126\cosh (4(x+y))-135\cosh (6(x+y))) 
\notag \\
&&+27(5\sinh (4(x+y))-4\sinh (2(x+y))))-81)  \label{29}
\end{eqnarray}

The other iterations contain large inputs and are not given. A computational
software program could help to calculate the other iterations up to any
order.

\begin{case}
RPSM: Now let us solve FZK equation by RPSM. The $k$-th residual function is
\end{case}

\begin{equation}
Res_{k}(x,y,t)=\frac{1}{8}(u_{k}^{2}(x,y,t)_{xyy}+\frac{1}{8}%
(u_{k}^{2}(x,y,t)_{xxx}+(u_{k}^{2}(x,y,t)_{x}  \label{30}
\end{equation}

To determine $f_{1}(x,y)$, we consider $k=1$%
\begin{equation}
Res_{1}(x,y,t)=D_{t}^{\alpha }u_{1}(x,y,t)+\frac{\partial }{\partial x}%
u_{1}^{2}(x,y,t)+\frac{1}{8}\frac{\partial ^{3}}{\partial x^{3}}%
u_{1}^{2}(x,y,t)+\frac{1}{8}\frac{\partial ^{3}}{\partial x\partial y^{2}}%
u_{1}^{2}(x,y,t)  \label{31}
\end{equation}%
since%
\begin{equation*}
u_{1}(x,y,t)=f(x,y)+f_{1}(x,y)\frac{t^{\alpha }}{\Gamma \lbrack 1+\alpha ]}
\end{equation*}%
thus%
\begin{eqnarray}
Res_{1}(x,y,t) &=&f_{1}(x,y)+2u_{1}(x,y,t)\frac{\partial }{\partial x}%
u_{1}(x,y,t)+\frac{1}{4}\frac{\partial ^{2}}{\partial y^{2}}u_{1}(x,y,t)%
\frac{\partial }{\partial x}u_{1}(x,y,t)  \notag \\
&&+\frac{1}{2}\frac{\partial }{\partial y}u_{1}(x,y,t)\frac{\partial ^{2}}{%
\partial x\partial y}u_{1}(x,y,t)+\frac{1}{4}u_{1}(x,y,t)\frac{\partial ^{3}%
}{\partial x\partial y^{2}}u_{1}(x,y,t)  \notag \\
&&+\frac{3}{4}\frac{\partial }{\partial x}u_{1}(x,y,t)\frac{\partial ^{2}}{%
\partial x^{2}}u_{1}(x,y,t)+\frac{1}{4}u_{1}(x,y,t)\frac{\partial ^{3}}{%
\partial x^{3}}u_{1}(x,y,t)  \label{32}
\end{eqnarray}%
from $D_{t}^{(k-1)\alpha }Res_{k}(x,y,0)=0,$ $0<\alpha \leq 1,$ $x\in I,$ $%
k=1,2,3,...$ and

\begin{eqnarray}
Res_{1}(x,y,0) &=&f_{1}(x,y)+2f(x,y)f(x,y)_{x}+\frac{1}{8}\left(
2f(x,y)_{yy}f(x,y)_{x}+4f(x,y)_{y}f(x,y)_{xy}+2f(x,y)(f(x,y)_{xyy}\right)  
\notag \\
&&+\frac{1}{8}6f(x,y)_{x}f(x,y)_{xx}+2f(x,y)f(x,y)_{xxx}  \label{33}
\end{eqnarray}%
for setting $Res_{1}(x,y,0)=0$ we obtain%
\begin{eqnarray}
f_{1}(x,y) &=&-2f(x,y)(f(x,y))_{x}-\frac{1}{4}(f(x,y))_{yy}(f(x,y))_{x}-%
\frac{1}{2}(f(x,y))_{y}(f(x,y))_{xy}-\frac{1}{4}f(x,y)(f(x,y))_{xyy}  \notag
\\
&&-\frac{3}{4}(f(x,y))_{x}(f(x,y))_{xx}-\frac{1}{4}f(x,y)(f(x,y))_{xxx}
\label{34}
\end{eqnarray}

Therefore, the first RPS aproximate solution is%
\begin{eqnarray}
u_{1}(x,y,t) &=&f(x,y)+\frac{1}{4\Gamma \lbrack 1+\alpha ]}t^{\alpha
}(-8f(x,y)f(x,y)_{x}-f(x,y)_{yy}f(x,y)_{x}-2f(x,y)_{y}f(x,y)_{xy}  \notag \\
&&-f(x,y)f(x,y)_{xyy}-3f(x,y)_{x}f(x,y)_{xx}-f(x,y)f(x,y)_{xxx})  \label{35}
\end{eqnarray}%
To obtain $f_{2}(x,y)$ substituting the second truncated series%
\begin{equation*}
u_{2}(x,y,t)=f(x,y)+f_{1}(x,y)\frac{t^{\alpha }}{\Gamma \lbrack 1+\alpha ]}%
+f_{2}(x,y)\frac{t^{2\alpha }}{\Gamma \lbrack 1+2\alpha ]}
\end{equation*}%
into the second residual function $Res_{2}(x,y,t)$%
\begin{equation}
Res_{2}(x,y,t)=D_{t}^{\alpha }u_{2}(x,y,t)+\frac{\partial }{\partial x}%
u_{2}^{2}(x,y,t)+\frac{1}{8}\frac{\partial ^{3}}{\partial x^{3}}%
u_{2}^{2}(x,y,t)+\frac{1}{8}\frac{\partial ^{3}}{\partial x\partial y^{2}}%
u_{2}^{2}(x,y,t)  \label{36}
\end{equation}%
Applying $D_{t}^{\alpha }$ on both sides and solving the equation $%
D_{t}^{\alpha }Res_{2}(x,y,0)=0$ gives

\begin{eqnarray}
f_{2}(x,y) &=&\frac{1}{4}%
(-8f_{1}(x,y)f(x,y)_{x}-f_{1}(x,y)_{yy}f(x,y)_{x}-8f(x,y)f_{1}(x,y)_{x}-f(x,y)_{yy}f_{1}(x,y)_{x}
\notag \\
&&-2f_{1}(x,y)_{y}f(x,y)_{xy}-2f(x,y)_{y}f_{1}(x,y)_{xy}-f_{1}(x,y)f(x,y)_{xyy}
\notag \\
&&-f(x,y)f_{1}(x,y)_{xyy}-3f_{1}(x,y)_{x}f(x,y)_{xx}-3f(x,y)_{x}f_{1}(x,y)_{xx}
\notag \\
&&-f_{1}(x,y)f(x,y)_{xxx}-f(x,y)f_{1}(x,y)_{xxx})  \label{37}
\end{eqnarray}%
so the first RPS aproximate solution is

\begin{eqnarray}
u_{2}(x,y,t) &=&f(x,y)+f_{1}(x,y)\frac{t^{\alpha }}{\Gamma (\alpha +1)}+%
\frac{1}{4\Gamma (2\alpha +1)}t^{2\alpha
}(-8f_{1}(x,y)f(x,y)_{x}-f_{1}(x,y)_{yy}f(x,y)_{x}  \notag \\
&&-3f_{1}(x,y)_{xx}f(x,y)_{x}-f(x,y)_{yy}f_{1}(x,y)_{x}-2f_{1}(x,y)_{y}f(x,y)_{xy}-2f(x,y)_{y}f_{1}(x,y)_{xy}
\notag \\
&&-f_{1}(x,y)f(x,y)_{xyy}-3f_{1}(x,y)_{x}f(x,y)_{xx}-f_{1}(x,y)f(x,y)_{xxx}-8f(x,y)f_{1}(x,y)_{x}
\notag \\
&&-f(x,y)f_{1}(x,y)_{xyy}-f(x,y)f_{1}(x,y)_{xxx})  \label{38}
\end{eqnarray}%
for the initial condition $f(x,y)=u\left( x,y,0\right) =\frac{4}{3}\rho
\sinh ^{2}(x+y)$. Following this manner the third iteration result, $%
u_{3}(x,y,t)$ is calculated. In Table 1, the third order approximate PIA and
RPSM results are compared numerically. Figure 1, Figure 2 and Figure 3 prove
that PIA and RPSM both give remarkably approximate results. We claim that
the higher iterations would give closer results. \bigskip 

\begin{table}[th]
\caption{Numerical results of Example 3.1. for different $u$ values when $%
\protect\alpha =1$ }
\begin{center}
{\scriptsize 
\begin{tabular}{cccccccccccc}
\hline
$x$ & $y$ & $t$ & \multicolumn{2}{c}{$\alpha =0.67$} & \multicolumn{2}{c}{$%
\alpha =0.75$} & \multicolumn{5}{c}{$\alpha =1.00$} \\ \hline
&  &  & $PIA$ & $RPSM$ & $PIA$ & $RPSM$ & $PIA$ & $RPSM$ & $Exact$ & $PIA$ $%
Error$ & $RPSM$ $Error$ \\ 
$0.1$ & $0.1$ & $0.2$ & 5.31854E-5 & 5.31244E-5 & 5.32747E-5 & 5.32479E-5 & 
5.35536E-5 & 5.35536E-5 & 5.39388E-5 & 3.85217E-7 & 3.85217E-7 \\ 
&  & $0.3$ & 5.28631E-5 & 5.28410E-5 & 5.29757E-5 & 5.29675E-5 & 5.33082E-5
& 5.33082E-5 & 5.38841E-5 & 5.75911E-7 & 5.75912E-7 \\ 
&  & $0.4$ & 5.25777E-5 & 5.25897E-5 & 5.27039E-5 & 5.27119E-5 & 5.30641E-5
& 5.30641E-5 & 5.38294E-5 & 7.65350E-7 & 7.65352E-7 \\ 
$0.6$ & $0.6$ & $0.2$ & 2.95493E-3 & 2.95185E-3 & 2.96356E-3 & 2.96251E-3 & 
2.98987E-3 & 2.98987E-3 & 3.03651E-3 & 4.66337E-5 & 4.66389E-5 \\ 
&  & $0.3$ & 2.92662E-3 & 2.92709E-3 & 2.93717E-3 & 2.93780E-3 & 2.96717E-3
& 2.96715E-3 & 3.03578E-3 & 6.86056E-5 & 6.86314E-5 \\ 
&  & $0.4$ & 2.90307E-3 & 2.90522E-3 & 2.91448E-3 & 2.91561E-3 & 2.94523E-3
& 2.94515E-3 & 3.03505E-3 & 8.98243E-5 & 8.99046E-5 \\ 
$0.9$ & $0.9$ & $0.2$ & 1.06822E-2 & 1.05506E-2 & 1.07716E-2 & 1.07143E-2 & 
1.10248E-2 & 1.10227E-2 & 1.15369E-2 & 5.12131E-4 & 5.14241E-4 \\ 
&  & $0.3$ & 1.04487E-2 & 1.01199E-2 & 1.05488E-2 & 1.03695E-2 & 1.07964E-2
& 1.07861E-2 & 1.15345E-2 & 7.38186E-4 & 7.48450E-4 \\ 
&  & $0.4$ & 1.02777E-2 & 9.60606E-3 & 1.03736E-2 & 9.96743E-3 & 1.05742E-2
& 1.05429E-2 & 1.15321E-2 & 9.57942E-4 & 9.89139E-4 \\ \hline
\end{tabular}
}
\end{center}
\end{table}

\begin{figure}[tbp]
\centering
\includegraphics[width=3.50in]{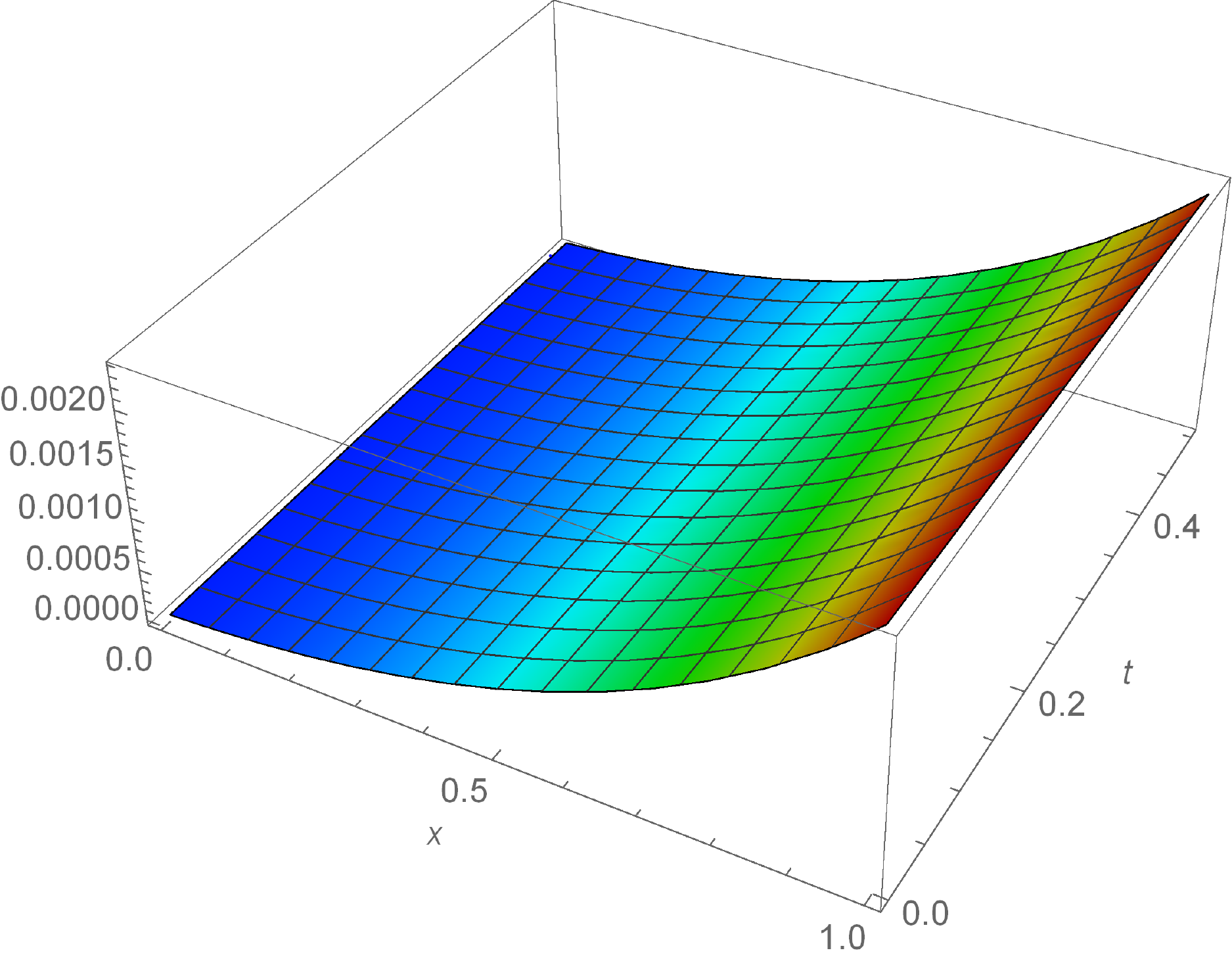}
\caption{PIA solution of fractional FZK Equation}
\end{figure}

\begin{figure}[tbp]
\centering
\includegraphics[width=3.50in]{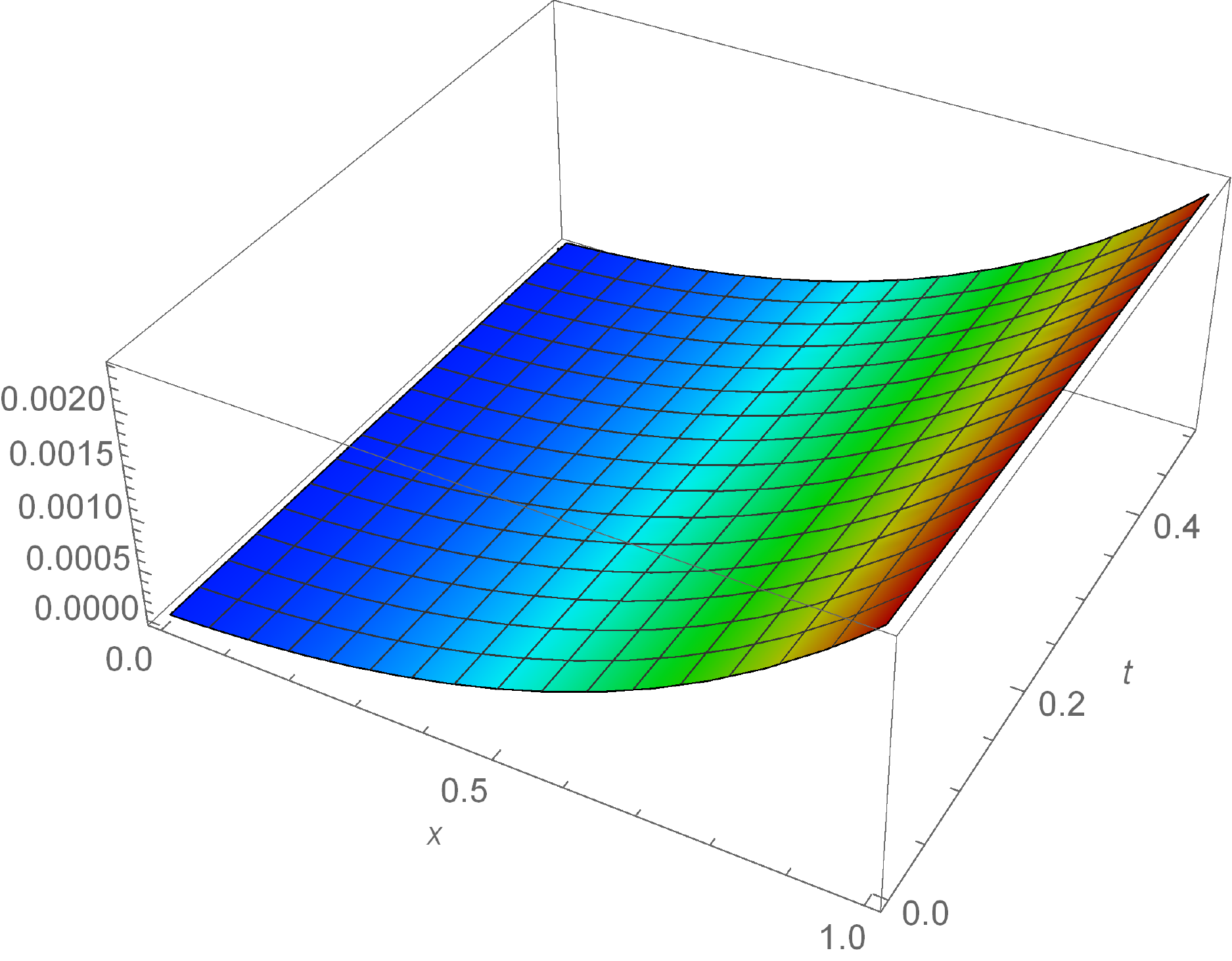}
\caption{RPSM solution of fractional FZK Equation}
\end{figure}

\begin{figure}[tbp]
\centering
\includegraphics[width=3.50in]{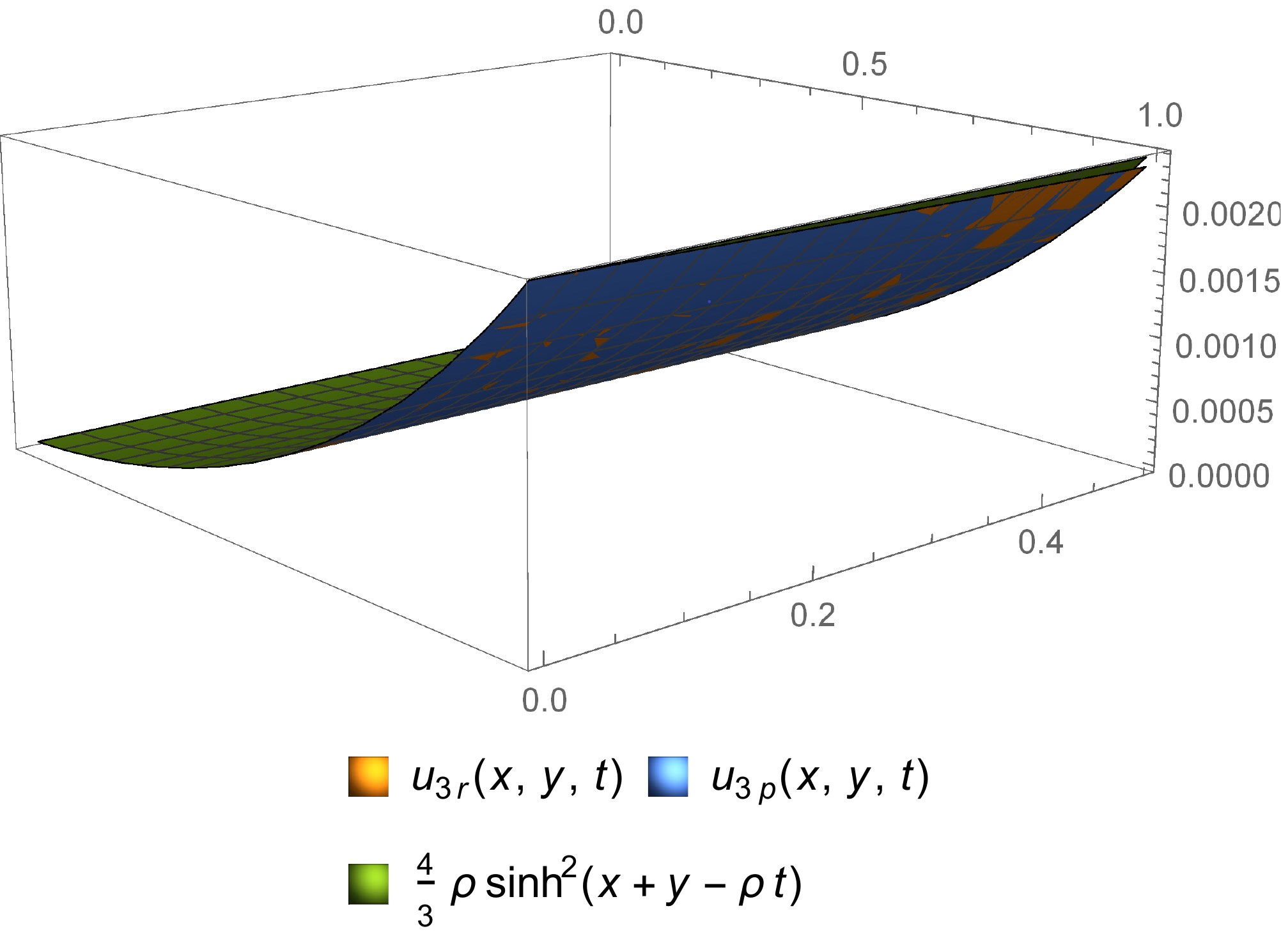}
\caption{Comparison of the PIA, RPSM and exact solutions of fractional FZK
Equation}
\end{figure}

\section{Conclusion}

\qquad In this study, perturbation-iteration algorithm and residual power
series method was introduced for time-fractional Zakharov-Kuznetsov partial
differential equation. It is clear that these methods are very simple and
reliable techniques and producing highly approximate results. We expect that
these methods could used to calculate the approximate solutions of other
types of fractional differential equations.

\end{document}